\documentclass[11pt]{article}

\usepackage{amsmath,latexsym,epsfig}
\usepackage{amsfonts,amssymb,amscd}

\title{\bf
Transitive factorisations in the symmetric group, and
combinatorial aspects of singularity theory
}

\author{I.P.Goulden\thanks{Dept. of Combinatorics and Optimization,
University of Waterloo, Waterloo, Ontario, Canada} and
D.M.Jackson\thanks{Dept. of Combinatorics and Optimization,
University of Waterloo, Waterloo, Ontario, Canada}
}
\date{ March 1999}

\begin{document}
\maketitle

 \newtheorem{theorem}{Theorem}[section]
 \newtheorem{proposition}[theorem]{Proposition}
 \newtheorem{definition}[theorem]{Definition}
 \newtheorem{axiom}[theorem]{Axiom}
 \newtheorem{lemma}[theorem]{Lemma}
 \newtheorem{corollary}[theorem]{Corollary}
 \newtheorem{remark}[theorem]{Remark}
\newtheorem{example}[theorem]{Example}
 \newtheorem{conjecture}[theorem]{Conjecture}

\def\cA{{\cal{A}}}
\def\cB{{\cal{B}}}
\def\cC{{\cal{C}}}
\def\cD{{\cal{D}}}
\def\cN{{\cal{D}}}
\def\cU{{\cal{U}}}
\def\cV{{\cal{V}}}

\def\bfx{{\rm\bf x}}
\def\bfi{{\rm\bf i}}

\def\frkt{\,{\mathfrak t}\,}

\def\sv{{\sf{v}}}
\def\sw{{\sf{w}}}

\def\proof{{\rm\bf Proof:\quad}}
\def\symgp{{\mathfrak S}}
\def\aut{{\rm{aut\,}}}
\def\bfone{{\bf 1}}

\def\qed{{\hfill{\large $\Box$}}}

\def\tsh{{\textstyle{\frac12}}}
\def\tsonethr{{\textstyle{\frac13}}}



\def\av#1{\left\langle#1\right\rangle_{\cW_N}}
\def\efM{e^{-\tsqu\tr\bfM^2}}
\def\etM{e^{-\tsh\tr\bfM^2}}
\def\ip#1#2{\left\langle #1,#2\right\rangle}
\def\avv#1{\left\langle#1\right\rangle_{\cV_N}}

\def\atp#1#2{\stackrel{\scriptstyle{#1}}{\scriptstyle{#2}}}

\begin{abstract}
We consider the determination of the number $c_k(\alpha)$ of ordered factorisations of an
arbitrary permutation on $n$ symbols, with cycle distribution $\alpha,$  
into $k$-cycles such that the factorisations have  minimal length and such that the 
 group generated by the factors acts transitively on the $n$ symbols.
The case $k=2$ corresponds to the
celebrated result of Hurwitz on the number of topologically distinct
holomorphic functions on the 2-sphere
that preserve a given number of elementary branch point singularities. 
In this case the monodromy group is the full symmetric group.
For  $k=3,$ the monodromy group is the  alternating group,
and this is another case that, in principle, is of considerable interest.

We conjecture an explicit form, for arbitrary $k,$ for the generating series for 
$c_k(\alpha),$ and prove that it holds for factorisations of permutations
with one, two and three cycles ($\alpha$ is a partition with at most three parts).
The generating series is naturally expressed in terms of the  symmetric functions 
dual to the those introduced by Macdonald for
the ``top'' connection coefficients in the class algebra of the symmetric group.

Our approach is  to determine a differential equation for the generating series
from a combinatorial analysis of the creation and annihilation of cycles
in products under the minimality condition.
\end{abstract}


\section{Introduction}\label{Si}
\subsection{Background}
This paper has two goals.
The first is to provide some general techniques to assist in the solution of
the type of permutation factorisation questions, with  transitivity 
and  minimality conditions, that originate in
the classical study of holomorphic mappings and branched coverings
of Riemann surfaces. Thus, we are concerned with certain {\em combinatorial}
questions that are encountered in aspects of {\em singularity theory}.
The appearance of such questions has long been recognized, and the reader
is directed to Arnold~\cite{arn}, for example, for further instances.

Very briefly, the classical construction concerns rational
mappings from a Riemann surface to the sphere.
Let $\alpha$ be the partition formed by the orders of the
poles of this  mapping. The poles are mapped to the point
at infinity.
Each  factor in an ordered factorisation is associated with a distinguished 
branch point, and it specifies the sheet transitions imposed
in a closed tour of the branch point, starting from an arbitrarily chosen
base point on the codomain of the mapping. 
In the generic case, the sheet transitions are
transposition ($2$-cycles).
The concatenation of the tours for each branch point,
from the same base point, in
the designated order, gives a sheet transition that is the product
of the sheet transitions for each branch point.
But this sheet transition is a permutation with $\alpha$ as its
cycle-type.
The transitivity condition ensures that the ramified covering
is connected, so the resulting Riemann surface is a ramified covering of a
sphere. The minimality condition ensures that the covering surface is 
a sphere also. 
The monodromy group is the group freely generated by the 
sheet transitions.

The particular class of permutation factorisation questions that
we shall consider in this paper involve as factors only $k$-cycles,
for some fixed, but arbitrary, value of $k$. The results that
we are able to obtain are thus extensions
of Hurwitz's~\cite{Hur} result with transpositions as factors, which
arose in the singularity theory context described above.

The second goal is to investigate the possibility of determining
 analogues of Macdonald's ``top'' symmetric functions that will
be appropriate for accommodating the transitivity condition. (It will
be recalled that Macdonald's top symmetric functions are associated
in a fundamental way with minimal ordered factorisations.) The possibility
of this connection arises from the fact that there is a striking common
element between the results of this paper on transitive minimal ordered
factorisations, and Macdonald's symmetric functions. This common element
is the functional equation

\begin{eqnarray}\label{equA}
w=xe^{w^{k-1}},
\end{eqnarray}
that arises in both settings when $k$-cycles are factors,
for apparently different reasons. The nature of this
possible connection is explored more fully in Section 1.5.

We will refer to these two contexts again, as the ideas in this paper
are developed. However, for the most part we now regard ordered
factorisations as discrete structures and we treat them by combinatorial
techniques. Throughout, we work in the appropriate ring of formal
power series. Thus, for example, the functional equation~(\ref{equA}) has
a unique solution for formal power series in $x$.
Although we have not completely attained the two goals, we have provided a substantial
amount of methodology for the first, and concrete evidence for the second.
We hope that the results are  substantial enough to provoke others
to explore further.

\subsection{Minimal ordered factorisations}

Let $\kappa(\pi)$ denote the number of
cycles in $\pi\in\symgp_n$. There is an obvious restriction
on $\kappa(\pi)$ under permutation multiplication. 

\begin{proposition}\label{Ppipi}
Let $\pi, \pi^\prime\in\symgp_n.$ Then
$$(n-\kappa(\pi))+(n-\kappa(\pi^\prime)) \ge (n-\kappa(\pi\pi^\prime)).$$
\end{proposition}

If $(\sigma_1,\ldots,\sigma_j)\in\symgp_n^j$ and $\sigma_1\cdots\sigma_j=\pi,$
then $(\sigma_1,\ldots,\sigma_j)$ is called an {\em ordered factorisation} of $\pi.$
Immediately from Proposition~\ref{Ppipi}, we obtain the inequality

\begin{equation}\label{bound}
\sum_{i=1}^j (n-\kappa(\sigma_i ))\geq n-\kappa(\pi ).
\end{equation}
In the case of equality, we call $(\sigma_1,\ldots,\sigma_j)\in\symgp_n^j$ a
{\em minimal} ordered factorisation of $\pi$.

Such factorisations have an elegant theory
and many enumerative applications (see, for example,
Goulden and Jackson~\cite{GJ}), including permissible commutation
of adjacent factors.
In particular, \cite{GJ} contains an explicit construction for a set
of symmetric functions (Macdonald's top symmetric functions)
that we shall return to in Section 1.5 of the Introduction.
Now we turn to the topic of the present paper.


\subsection{Minimal transitive ordered factorisations}

We write $\alpha\vdash n$ to indicate that $\alpha$ is a partition of $n,$ and
$\cC_\alpha$ for the conjugacy
class of $\symgp_n$ indexed by $\alpha.$ Let $l(\alpha)$ denote the
number of parts in $\alpha$.
If  $\pi\in\cC_\alpha$ then $\kappa(\pi)=l(\alpha).$
An ordered factorisation $(\sigma_1,\ldots,\sigma_j)$ is said to be
{\em transitive} if the subgroup of $\symgp_n$ generated by the factors
acts transitively on $\{1,\ldots,n\}.$ The case where each of the factors
is in $\cC_{[k,1^{n-k}]},$ and is therefore a pure $k$-cycle, is of
particular interest.
A transitive ordered factorisation of $\pi\in\cC_\alpha$ with the minimal
choice of $j$ consistent with the other conditions
is said to be {\em minimal}. In this case, $j=\mu_k(\alpha),$ 
where
$$\mu_k(\alpha)=\frac{n+l(\alpha)-2}{k-1},$$
as we shall prove in Proposition~\ref{Pmu}.
For example, when $k=3$,

\begin{eqnarray}\label{threxamp}
(247)(586)(479)(136)(235)=(1386)(254)(79),
\end{eqnarray}
and $((247),(568),(479),(136),(235))$ is a minimal transitive ordered
factorisation of the permutation 
$(1386)(254)(79)$, into $3$-cycles with $5$ factors
(minimality holds in this example since $\mu_3([4,3,2])=5$).

Such factorisations are encountered in a number of contexts.
These include, for example, the topological classification of polynomials  of given degree and a 
given number of critical values, and the moduli space of covers of the Riemann sphere and properties
of the Hurwitz monodromy group, and applications to
mathematical physics~\cite{crt}. 
The reader is directed to~\cite{fried,Zvon,KZ} for
further background information.

The number
of minimal transitive ordered factorisations of an arbitrary
but fixed $\pi\in\cC_\alpha$ is denoted by $c_k(\alpha).$
Hurwitz~\cite{Hur} determined
$c_2(\alpha)$, as a consequence of his study of holomorphic
mappings on the sphere (see also Strehl~\cite{str97}, for the proof of
an identity that completes Hurwitz's treatment). He showed that
\begin{eqnarray}\label{eHres}
c_2(\alpha)=n^{l(\alpha)-3}(n+l(\alpha )-2)!\prod_{j=1}^{l(\alpha)}
\frac{\alpha_j^{\alpha_j}}{(\alpha_j-1)!}.
\end{eqnarray}
A shorter and self-contained proof of this result has been given
by Goulden and Jackson~\cite{GJ97}.
The special case $c_2([1^n])$
was derived independently by Crescimanno and Taylor~\cite{crt}.
For related work, in the language of singularity theory,
see~\cite{ssvain}.

The case $k=3$ is also of considerable interest, for the subgroup
generated is the alternating group.

\subsection{The results and a  conjecture}
The main conjecture of the paper
concerns the form of the generating series for the
$c_k(\alpha)$.
Let $u,z,p_1,p_2,\ldots$ be
indeterminates and let $p_\alpha=p_{\alpha_1}p_{\alpha_2}\ldots$. Then 
\begin{eqnarray*}
F_k^{(m)}(u,z;p_1,p_2,\ldots) &=&\sum_{\atp{n\ge1}{k-1\vert n+m-2}}
\sum_{\atp{\alpha\vdash n}{ l(\alpha)=m}}
c_k(\alpha)\,\vert\cC_\alpha\vert\, p_\alpha
\frac{u^{\mu_k(\alpha)}}{\mu_k(\alpha)!}
\frac{z^n}{n!}.
\end{eqnarray*}
The series $F_k^{(m)}$ is a formal power series in $z$ with coefficients that are 
polynomial in $u,p_1,p_2,\ldots,$ and we will be working in this ring.

It is more convenient to work with a symmetrised form
of the generating series, defined in terms of the following operator
 $\psi_m.$ 
If $\alpha$ is a partition with $m$ parts, then
\begin{equation}\label{psidef}
\psi_m\left(p_\alpha u^iz^j\right)=\sum_{\sigma\in\symgp_m}
x_1^{\alpha_{\sigma(1)}}\cdots x_m^{\alpha_{\sigma(m)}}.
\end{equation}
Now define
$$P_k^{(m)}(x_1,\ldots,x_m)=\psi_m(F_k^{(m)}).$$

In the main conjecture that follows, we let  $w_i=w(x_i)$ for $i\ge1$,
and $w(x)$ is the
unique power series solution of the functional equation given
in~(\ref{equA}).

\begin{conjecture}\label{CConj}
For $m\ge1,$
$$
\left( \sum_{i=1}^m x_i\frac{\partial}{\partial x_i}\right)^{3-m}
P_k^{(m)}(x_1,\ldots,x_m) =
S_k^{(m)}(w_1,\ldots,w_m)\,
\prod_{i=1}^m x_i\frac{d w_i}{d x_i},
$$
where $S^{(m)}_k(w_1,\ldots,w_m)$ is  a symmetric polynomial in $w_1,\ldots,w_m.$
\end{conjecture}

The conjectured form for the series $P_k^{(m)}$ therefore involves rational
expressions in $w_1,\ldots,w_m.$  To see this,  
differentiate~(\ref{equA}) with respect to $x$, 
to obtain the  rational form
\begin{eqnarray}\label{equAA}
x \frac{dw}{dx} = \frac{w}{1-(k-1)w^{k-1}}.
\end{eqnarray}
Note that the dependence on $k$ rests in the coefficients of the symmetric
polynomial (which we conjecture to be polynomials in $k$), but more
deeply in the functional equation~(\ref{equA}).
The explicit formal power series for $w$ is actually straightforward, and
obtained immediately
by Lagrange's Theorem, yielding
\begin{eqnarray}\label{Ewx}
w(x) = \sum_{m\ge0}\frac{(1+(k-1)m)^{m-1}}{m!}x^{1+(k-1)m}.
\end{eqnarray}

\medskip
In this paper, we are able to determine explicitly $P_k^{(m)}$ for
the cases $m=1,2,3$. These are all of a form that satisfies the above
conjecture. The resulting expressions for $S_k^{(m)}$ in these cases
are stated below.
Let $V(w_1,\ldots,w_m)$ denote the Vandermonde determinant in $w_1,\ldots,w_m.$

\begin{theorem}\label{Tonept}
$S^{(1)}_k(w_1)=1.$
\end{theorem}


\begin{theorem}\label{Ttwopt}
$S^{(2)}_k(w_1,w_2) =(w_1^{k-1}-w_2^{k-1})^2/V(w_1,w_2)^2.$
\end{theorem}

\begin{theorem}\label{Tthreept}
$S^{(3)}_k(w_1,w_2,w_3) =G^2/V(w_1,w_2,w_3)^2,$ where
\begin{eqnarray*}
G&=& w_1\left(1-(k-1)w_1^{k-1}\right) (w_3^{k-1}-w_2^{k-1}) 
+ w_2\left(1-(k-1)w_2^{k-1}\right) (w_1^{k-1}-w_3^{k-1}) \\
&\mbox{}&+ w_3\left(1-(k-1)w_3^{k-1}\right) (w_2^{k-1}-w_1^{k-1}).
\end{eqnarray*}
\end{theorem}

\medskip
The proofs of these results are given in Section~4 of
the paper. The method is to solve a partial differential equation
for $P_k^{(m)}$ that is obtained in Section~3. This
equation is itself deduced by symmetrising a partial differential
equation for $F_k^{(m)}$ that is obtained in Section~2. The latter
is determined by a combinatorial analysis of minimal permutation
multiplication.

The determination of further cases, at present, seems
to be intractable, as we discuss in Section~5.
The forms obtained above in the first three cases
are remarkably simple, although it has not been
possible to conjecture a sufficiently precise general form based on this evidence.
Although $S^{(1)}_k,$ by default, $S^{(2)}_k$ and $S^{(3)}_k$   are perfect squares,
we do not believe that this holds in general. 

Note that
$S^{(m)}_2$ does not restrict to $S^{(m-1)}_2$ through $w_m=0,$ in the
cases $m=2$ and $m=3$.
Also note that if we substitute $k=2$ in Theorems~\ref{Ttwopt} and~\ref{Tthreept}
above, then we immediately obtain $S^{(2)}_2=S^{(3)}_2=1$. In the following
result, we demonstrate that this is true when $k=2$ for arbitrary choice
of $m$, as a direct consequence of Hurwitz's result.

\begin{lemma}\label{LSm2}
$S^{(m)}_2(w_1,\ldots,w_m)=1$ for $m\ge1.$
\end{lemma}
\proof From~(\ref{eHres}),
\begin{eqnarray*}
\left( \sum_{i=1}^m x_i\frac{\partial}{\partial x_i}\right)^{3-m}
P_2^{(m)}(x_1,\ldots,x_m) &=&
\sum_{n\ge1}\sum_{\atp{\alpha\vdash n}{l(\alpha)=m}}
\frac{\vert\cC_\alpha\vert}{n!}
\left(\prod_{j=1}^m\frac{\alpha_j^{\alpha_j+1}}{\alpha_j!}\right)
\sum_{\sigma\in\symgp_m}
x_1^{\alpha_{\sigma(1)}}\cdots x_m^{\alpha_{\sigma(m)}} \\
&=& \frac{1}{m!}\sum_{\alpha_1,\ldots,\alpha_m\ge1}
\left(\prod_{j=1}^m\frac{\alpha_j^{\alpha_j}}{\alpha_j!}\right)
\sum_{\sigma\in\symgp_m}
x_1^{\alpha_{\sigma(1)}}\cdots x_m^{\alpha_{\sigma(m)}} \\
&=& \prod_{j=1}^m x_j\frac{dw_j}{dx_j}.
\end{eqnarray*}
The result now follows. \qed

\medskip
We note that, in the case of transpositions, together with
Vainshtein~\cite{gjvain}, we have recently been
able to obtain similar results in the case where there are two more than the
minimal number of factors. These correspond to holomorphic mappings
from the torus.

\subsection{Symmetric functions and minimal ordered factorisations}

In~\cite{GJ}(see also~\cite{MAC}) an explicit construction is given
for symmetric functions $u_\lambda,$ indexed by $\lambda\vdash n,$
that are closely related to minimal ordered factorisations in
the symmetric group (note that the term ``top'' was used for
such factorisations in that paper; these are Macdonald's
{\em top symmetric functions}). In particular,
the number of minimal ordered
factorisations $(\sigma_1, \ldots ,\sigma_j)$ of $\pi$, where
$\sigma_i\in\cC_{\beta_i}, i=1,\ldots,j$, and for each $\pi\in\cC_\lambda$,
is given by
\begin{equation}\label{symmu}
[u_{\lambda-1}]\,u_{\beta_1-1}\cdots u_{\beta_j-1},
\end{equation}
where $\beta_i-1$ is the partition obtained by subtracting one from each
part of $\beta_i.$ Properties that can be developed for $u_\lambda$
then facilitate the determination of this number. Several examples
of their use in enumerative questions are given in~\cite{GJ},
together with the enumeration of minimal ordered factorisations
up to permissible commutation of adjacent factors.
 
We now recall the algebraic construction for the symmetric functions
$u_\lambda,$ where $\lambda\vdash n.$
Let $H(t;\bfx)$ be the generating series for the complete symmetric
functions $h_k(\bfx)$ of degree $k$ in $\bfx=(x_1,x_2,\ldots).$
Then the functional equation $s=t\,H(t;\bfx)$ has a unique solution
$t\equiv t(s,\bfx)$  given by $t=s\,H^\star (s;\bfx)$ where
$H^\star(s;\bfx)=\sum_{j\ge0}s^j h^\star_j(\bfx),$ and $h^\star_j(\bfx)$
is a symmetric function in $\bfx$ of total degree $j.$ Let
$h_\lambda^\star=h_{\lambda_1}^\star h_{\lambda_2}^\star \cdots.$
Then  $\{u_\lambda\}$ is defined to be the basis for the symmetric
function ring that is dual to the  basis $\{h_\lambda^\star\}$
with respect to the inner product for which the monomial and
complete symmetric functions are dual (see, e.g.
Macdonald~\cite{MAC}, for a complete treatment of the required background
material).

Thus, for minimal ordered factorisations in which all
factors are $k$-cycles, then in equation~(\ref{symmu}), we
have $u_{\beta_i -1} =u_{k-1}$ for all $i=1,\ldots ,j$.
But, as is shown in~\cite{GJ}, $u_{k-1}=-p_{k-1}$,
so for minimal ordered factorisations in which all
factors are $k$-cycles, we can restrict attention to a
symmetric function algebra in which $p_i=0$ if $i\neq k-1$.
In this case, we have
$$ s=t H(t;\bfx)= \exp\left(\sum_{m\geq 1} \frac{p_m}{m}t^m\right)=t \exp
\left(\frac{-p_{k-1}}{k-1}t^{k-1}\right).
$$
Thus, if $z$ is substituted for $\frac{p_{k-1}}{k-1}$, in this equation,
we obtain
$$t=s e^{zt^{k-1}}.$$
But this is precisely the functional equation~(\ref{equA}), whose
solution features so centrally
in our results
for the transitive case above.

We conclude from this that there must be an important relationship between the
transitive case of minimal ordered factorisations for which we have obtained
partial results in this paper, and minimal ordered factorisations
themselves, that have such an elegant theory based on
symmetric functions. Although we have been unable to find
a direct link between these two classes, we hope that the results of
this paper will provide a good starting point for such a direct link,
and a similarly elegant
theory for the transitive case.



\section{The partial differential equation}\label{Sconst}
In this section we determine a partial differential equation for the
generating series

\begin{eqnarray*}
\Phi^{(k)} &=& \sum_{m\ge1}F_k^{(m)}
\end{eqnarray*}
by  a case analysis of the creation and annihilation of cycles
in products of permutations subject to the minimality condition.

We begin with a discussion of permutation multiplication. First,
we prove
the expression that has been given in Section 1.3 for $\mu_k(\pi)$.

\begin{proposition}\label{Pmu}
Let $\alpha\vdash n,$ and let $\pi\in\cC_\alpha.$ Then $\mu_k(\pi)=\mu_k(\alpha),$ where
$$\mu_k(\alpha)=\frac{n+l(\alpha)-2}{k-1}.$$
\end{proposition}
\proof Let $(\sigma_1,\ldots,\sigma_j)$ be  a minimal transitive ordered factorisation
of $\pi$ into $k$-cycles. Let $\pi^\prime$ and $\pi$ be in the same conjugacy class, so
$\pi^\prime=g^{-1}\pi g$ for some $g\in\symgp_n.$ Then
$(g^{-1}\sigma_1 g,\ldots,g^{-1}\sigma_j g)$  is a minimal transitive ordered factorisation
of $\pi^\prime$, so $\mu_k(\pi^\prime)=\mu_k(\pi),$ and we denote the common value
by $\mu_k(\alpha)$ where $\pi\in\cC_\alpha.$
Now each $k$-cycle in  $\symgp_k$ has a minimal transitive ordered factorisation
into $\mu_2([k])$ transpositions, so
$\mu_2(\alpha)=\mu_2([k])\,\mu_k(\alpha).$
But (Prop.~2.1,~\cite{GJ97}), $\mu_2(\alpha)=n+l(\alpha)-2,$ and the result
follows. \qed

\medskip
Next we give a combinatorial characterisation of minimal transitive ordered
factorisations.
The following lemma  characterises the relationship between $\sigma_1$ and
$\sigma_2\cdots\sigma_j$ for a minimal transitive ordered factorisation
$(\sigma_1,\ldots,\sigma_j)$ of $\pi\in\symgp_n$ into $k$-cycles.
Some terminology will be useful.
The multi-graph $\cD_{\sigma_1,\ldots,\sigma_j}$ has vertex-set
 $\{1,\ldots,n\}$, and edges consisting of the edges of the $k$-cycles in the
factorisation.
Let $\cV_1,\ldots,\cV_l$ be
the vertex-sets of the connected components
of $\cD_{\sigma_2,\ldots,\sigma_j}$,
so $\{\cV_1,\ldots,\cV_l\}$ is a
partition of $\{1,\ldots,n\}$ into nonempty subsets.
For $i=1,\ldots ,l$, let $\alpha_i$ consist of all $t\in\{ 2,\ldots ,j\}$ such
that all of the $k$ elements on $\sigma_t$ belong to $\cV_i$,
so $\{ \alpha_1,\ldots ,\alpha_l\}$ is a partition of $\{ 2,\ldots ,j\}$.
Suppose $\alpha_i = \{ \alpha_{i1},\ldots ,\alpha_{is_i}\}$,
with $\alpha_{i1}<\cdots <\alpha_{is_i}$, and
 $\sigma_{\alpha_{i1}}\cdots \sigma_{\alpha_{is_i}}=\pi_i$, for $i=1,\ldots,l$.
Then clearly, by construction, $(\sigma_{\alpha_{i1}},\cdots 
,\sigma_{\alpha_{is_i}})$ is a minimal transitive ordered factorisation
of $\pi_i$, for $i=1,\ldots ,l$, and we have 

\begin{equation}\label{pisig}
\pi = \sigma_1 \pi_1\cdots
\pi_l.
\end{equation}

For example, in the minimal transitive factorisation given in~(\ref{threxamp}),
we have $l=2$, with $\cV_1 = \{ 1,2,3,5,6,8\}$, and $\cV_2 =\{ 4,7,9\}$;
 $\alpha_1 =\{ 2,4,5\}$, and $\alpha_2 =\{ 3\}$; $\pi_1 =(1386)(25)$,
and $\pi_2 =(479)$.

For $\pi\in\symgp_n$ and $\cA\subseteq\{1,\ldots,n\},$ the $\cA$-{\em restriction}
of $\pi$ is the permutation on $\cA$ obtained by deleting the elements not in $\cA$
from the cycles of $\pi.$
For example, if $\pi = (1538)(27469)$ and $\cA = \{ 1,4,6,7,8\}$,
then the $\cA$- restriction of $\pi$ is $(18)(467)$.

\begin{lemma}\label{Lmtof}
Let $(\sigma_1,\ldots,\sigma_j)$ be a minimal transitive ordered
factorisation of $\pi\in\symgp_n$ into $k$-cycles, and let
 $\pi_1,\ldots ,\pi_l$ be constructed as above. Then
\begin{enumerate}
   \item  $\sigma_1$ has at least one element in
 common with each of $\pi_1,\ldots,\pi_l.$
   \item The elements of $\sigma_1$ in common with $\pi_i$ lie on a single cycle
         of $\pi_i,$ for $i=1,\ldots,l.$
   \item Let $\cU$ denote the $k$-subset of $\{1,\ldots,n\}$ consisting
         of the elements on the $k$-cycle $\sigma_1$. Let $\gamma$ denote
         the $\cU$-restriction of $\sigma_1$, and let $\tau$ denote
         the $\cU$-restriction of $\pi$. If $\rho=\gamma^{-1}\tau$, then
         $(k-\kappa(\tau))+(k-\kappa(\rho))=k-\kappa(\gamma)$, so
         $(\tau ,\rho^{-1})$ is a minimal ordered factorisation
         of $\gamma$.
\end{enumerate}
\end{lemma}
\proof
Since $(\sigma_1,\ldots,\sigma_j)$ is a transitive factorisation
of $\pi$, then $\cD_{\sigma_1,\ldots,\sigma_l}$ is connected.
Thus the single $k$-cycle in  $D_{\sigma_1}$ has
at least one vertex in each of the 
connected components of $\cD_{\sigma_2,\ldots,\sigma_l}$, and
this establishes part 1.

Now, from~(\ref{pisig}) and the fact that $(\sigma_{\alpha_{i1}},\cdots
,\sigma_{\alpha_{is_i}})$ is a minimal transitive ordered factorisation
of $\pi_i$, for $i=1,\ldots ,l$, we have

\begin{equation}\label{musum}
\mu(\pi)=1+\mu(\pi_1)+\cdots+\mu(\pi_l).
\end{equation}
But, from Proposition~\ref{Pmu}
$$\mu(\pi)=\frac{n+\kappa(\pi)-2}{k-1}\qquad\mbox{and}\qquad
\mu(\pi_i)=\frac{\vert\cV_i\vert+\kappa(\pi_i)-2}{k-1},$$
for $i=1,\ldots ,l$.
Thus, substituting these values for the $\mu 's$ into~(\ref{musum})
we obtain 

\begin{eqnarray*}
n+\kappa(\pi)-2 = k-1+\sum_{i=1}^l (\vert\cV_i\vert+\kappa(\pi_i)-2). 
\end{eqnarray*}
But $n=\sum_{i=1}^l \vert\cV_i\vert$, and substituting
this into the above gives
\begin{eqnarray}\label{se2}
\kappa(\pi)-\sum_{i=1}^l\kappa(\pi_i)=k+1-2l.
\end{eqnarray}
Now let $\rho_i$ be the $\cU$-restriction of $\pi_i$, for $i=1,\ldots ,l$,
so $\pi=\sigma_1\pi_1 \cdots\pi_l$ restricts down to $\tau=\gamma\rho$,
where $\rho=\rho_1 \cdots\rho_l$. We then have
$$\kappa(\pi)=\kappa(\tau)+\sum_{i=1}^l (\kappa(\pi_i)-\kappa(\rho_i))
 = \kappa(\tau)+\sum_{i=1}^l \kappa(\pi_i) - \kappa(\rho),$$
and together with~(\ref{se2}) this gives
\begin{equation}\label{diffkl}
\kappa(\tau)-\kappa(\rho)= \kappa(\pi)-\sum_{i=1}^l\kappa(\pi_i)=k+1-2l.
\end{equation}
On the other hand, since $\gamma, \rho$ and $\tau$ act on a $k$-set 
 and $\tau\rho^{-1} = \gamma$ we have from Proposition~\ref{Ppipi}  that  
$\left(k-\kappa(\tau)\right) + \left(k-\kappa(\rho^{-1})\right) \ge \left(k-\kappa(\gamma)\right).$
But $\kappa(\gamma)=1$ and $\kappa(\rho^{-1})=\kappa(\rho)$, so
$\kappa(\tau)+\kappa(\rho)\le k+1$, and in addition, from part 1 we have
$\kappa(\rho)\ge l.$ It follows that
$\kappa(\tau)-\kappa(\rho)\le k+1-2\kappa(\rho) \le k+1-2l.$
Combining this with~(\ref{diffkl}) gives $\kappa(\rho)=l$.
Together with part~1, this establishes part~2.

Part 3 follows immediately from
$\kappa(\rho)=l$, $\kappa(\gamma)=1$ and~(\ref{diffkl}). \qed

\medskip
We  now use this characterisation as a construction
for deriving a partial differential
equation for $\Phi^{(k)}$ with arbitrary $k$.
In the interests of succinctness, we suppress the occurrences  of $k$
in $\Phi^{(k)}$ and $P^{(m)}_k$.
From Lemma~\ref{Lmtof}(.3),
the terms in the equation
are in one-to-one correspondence with
minimal ordered factorisations of a $k$-cycle. These factorisations are
themselves in one-to-one correspondence with a particular class of trees,
as was shown in~\cite{GJcactus}, and described as follows:
Let $\cB^{(k)}$ be the set of all plane two-coloured (black, white) trees
with $k$ edges,
with the indices $i_1,\ldots,i_k$ assigned to different edges
in a canonical way. Let $\frkt$ be such a tree.
For $v\in\frkt$ let $\omega(v)$ be the sum of the indices of edges incident with $v.$
Let $\widehat{\frkt}$ denote the tree obtained from $\frkt$ by deleting monovalent white vertices.
Let $\aut(\widehat{\frkt})$ denote the automorphism group of $\widehat{\frkt}$ with the convention
that if $\widehat{\frkt}$ is an isolated black vertex, then $\aut(\widehat{\frkt})$ is the cyclic
group on $k$ symbols. 

\begin{theorem}\label{Tdiff1}
Let $\bfi=(i_1,\ldots,i_k)$ where $i_1,\dots,i_k\ge1.$ 
Then $\Phi$ satisfies the nonlinear, inhomogeneous partial differential equation
\begin{equation}\label{phidiffnew}
\sum_{\bfi\ge\bfone}
\sum_{\frkt\in\cB^{(k)}}
\frac{1}{\vert\aut(\widehat{\frkt})\vert}\,
\left(\prod_{\sv\in \cV_{black}(\frkt)} p_{\omega(\sv)} 
\prod_{\sw\in \cV_{white}(\frkt)}
\,
 \omega(\sw)\frac{\partial \Phi}{\partial p_{\omega(\sw)}}\right)
 =\frac{\partial\Phi}{\partial u},
\end{equation}
with the convention that empty sums are zero and empty products are equal to one.
\end{theorem}
\proof
From Lemma~\ref{Lmtof}(.3), $(\tau,\rho^{-1})$ is a minimal ordered factorisation
in $\symgp_k$ of the $k$-cycle $\gamma$. Thus, from~\cite{GJcactus} Theorem~2.1, 
$(\tau,\rho^{-1})$ uniquely encodes
an edge-rooted 2-coloured plane tree $\frkt$  with $k$ edges, such that the 
black vertex-degrees are given by the cycle-type of $\rho^{-1},$ and the white
vertex-degrees are given by the cycle-type of $\tau.$

We now observe that,
in the product $\gamma\rho,$ cycles with length equal to the degree of each
of the black vertices are annihilated, and combined to form cycles of length equal
to the degree of a white vertex.
This observation permits us to reconstruct the cycle distribution of
$\pi$ from $\sigma_1$ and the cycle distributions of $\pi_1,\ldots,\pi_l.$

The tree $\frkt$ can be regarded as the boundary of a polygon.
As the boundary of $\frkt$ is traversed, each edge is encountered twice, once in the direction
from its black vertex towards its white vertex, and once in the direction from its white vertex
towards its black vertex. The indexed symbols $i_1,\ldots,i_k$ are assigned
to the edges of $\frkt,$ starting from the root-edge, as each edge is encountered
in the direction from its black vertex towards its white vertex.
Moreover, $i_j$ is the number
of elements in $\{1,\ldots,n\}$ that separate two elements in $\sigma_1$ 
on cycles in $\pi_i.$ 

In this encoding the degree of a black vertex is the number of elements of $\sigma_1$
that are incident with $\pi_i$ and the number of black vertices is $\kappa(\rho).$
This indicates which cycles in $\pi_1,\ldots,\pi_l$ are annihilated in premultiplication
by $\sigma_1$ and which cycles are created. It is necessary only to keep track of
the lengths of these cycles.

The  contribution from cycles that are created is therefore
\begin{eqnarray*}
\prod_{\sv\in \cV_{black}(\frkt)} p_{\omega(\sv)}.
\end{eqnarray*}

The contribution from cycles that are annihilated is
\begin{eqnarray*}
\prod_{\sw\in \cV_{white}(\frkt)}
 \omega(\sw)\frac{\partial \Phi}{\partial p_{\omega(\sw)}}.
\end{eqnarray*}
To see this, select one of the cycles $\rho_i.$ Next select an element on it.
Then mark off the cycle into a number of contiguous segments equal to
the degree of the corresponding black vertex  in $\frkt.$
However, this overcounts by a factor of 
${1}/{\vert\aut(\widehat{\frkt})\vert}.$
Thus summing we have
\begin{eqnarray*}
\sum_{\bfi\ge\bfone}
\sum_{\frkt\in\cB^{(k)}}
\frac{1}{\vert\aut(\widehat{\frkt})\vert}\,
\left(\prod_{\sv\in \cV_{black}(\frkt)} p_{\omega(\sv)} \,
\prod_{\sw\in \cV_{white}(\frkt)}  \,
 \omega(\sw)\frac{\partial \Phi}{\partial p_{\omega(\sw)}}\right).
\end{eqnarray*}

But this is equal to the generating series for minimal transitive ordered factorisations
with the leftmost factor deleted.
But this is
${\partial\Phi}/{\partial u}.$
The result now follows. \qed

\medskip
Note that, if  $p_i$ is the power sum symmetric function of degree $k$ in an infinite
set of ground variables, then
$j {\partial}/{\partial p_{j}}=p^\star_{j},$ where $p^\star_{j}$ is
the adjoint of premultiplication by $p_{j}$ (see, e.g.,~\cite{MAC} for details).
The partial differential equation therefore can be rewritten in the following form,
that exhibits the symmetry between black and white vertices, as
\begin{eqnarray*}
\sum_{\bfi\ge\bfone}
\sum_{\frkt\in\cB^{(k)}}
\frac{1}{\vert\aut(\widehat{\frkt})\vert}\,
\left(\prod_{\sv\in \cV_{black}(\frkt)} p_{\omega(\sv)}\right)\,
\left(\prod_{\sw\in \cV_{white}(\frkt)}
  p^\star_{\omega(\sw)} \Phi\right)\,
 =\frac{\partial\Phi}{\partial u}.
\end{eqnarray*}


It will be useful to list explicitly the first few trees on
the left hand side of~(\ref{phidiffnew}) in the arbitrary case,
graded by the number of black vertices in $\widehat{\frkt}$,
to find the equations for the low order terms of $\Phi$, in the $p$'s. We consider
below all of the trees with at most three black vertices.

\medskip
\noindent{\bf First tree:} Let $\widehat{\frkt_1}$ be the tree in $\cB_k$ consisting of one black
vertex joined to $k$ white vertices. Assign $i_1,\ldots,i_k$ to the edges. Then
$\widehat{\frkt_1}$ is the tree consisting of an isolated black vertex.
By the convention on automorphisms, $\aut(\frkt_1)=k.$ 

\noindent{\bf Second tree:} Let $\widehat{\frkt_2}$ be the tree in $\cB_k$ consisting of 
a path $\widehat{\frkt_2}$ with two black vertices and one white vertex. 
Attach $i_{k-1}$ and $i_k$ to the two edges incident with the white vertex.
Now  join
$r$ white vertices  to one of the black vertices, and attach $i_1,\ldots, i_r$ to
the edges. Join $k-r-2$
white vertices  to the other black vertex, and attach labels $i_{r+1},\ldots, i_{k-2}$
to the edges. The resulting tree $\frkt_2$ therefore
has $k$ edges, and it is readily seen that
$\widehat{\frkt_2}$ is the tree obtained from $\frkt_2$ by removing monovalent white vertices.  
Moreover,  $\vert\aut(\widehat{\frkt_2})\vert=2.$

For the trees with three black vertices,
we give only the tree from which monovalent white vertices have been
removed. The summation variables $i_1,\ldots,i_{k}$ are attached to edges
in the way described in the previous cases.

\noindent{\bf Third tree:} Let $\widehat{\frkt_3}$ be the tree consisting of
one white vertex to which three black vertices are joined.
Then $\vert\aut(\widehat{\frkt_3})\vert=3.$ Note that, in this case,
the path separates into two sets the additional white vertices that are joined
to the  black vertex in the middle of the path.

\noindent{\bf Fourth tree:} Let $\widehat{\frkt_4}$ be the path consisting of
three black vertex and two white vertices.
Then $\vert\aut(\widehat{\frkt_4})\vert=2.$ 

\medskip
The partial  differential  equations for minimal transitive ordered factorisations
into 2-cycles, and into 3-cycles, can be written down explicitly from the
terms that have been given.
Let $\Phi_j\equiv j\partial \Phi/ \partial p_j$, for$j\geq 1$. 

When $k=2$  the only trees with two edges correspond
to $\widehat{\frkt_1}$ and$\widehat{\frkt_2}$, so
in this case~(\ref{phidiffnew}) becomes
\begin{equation}\label{twocase}
\tsh\sum_{i_1,i_2\ge1}
 \left(\Phi_{i_1}^{(2)}\Phi_{i_2}^{(2)}p_{i_1+i_2}+\Phi_{i_1+i_2}^{(2)} p_{i_1}p_{i_2}\right)=
\frac{\partial \Phi^{(2)}}{\partial u}.
\end{equation}
This is the equation given in~\cite{GJ97}, where we demonstrated that
a series conjectured from numerical computations satisfied the  equation uniquely.

When $k=3$  the only trees with three edges correspond
to $\widehat{\frkt_1}$, $\widehat{\frkt_2}$ and $\widehat{\frkt_3}$, so
\begin{equation}\label{thrpde}
\sum_{i_1,i_2,i_3\ge1}
 \left(  \tsonethr\Phi_{i_1}^{(3)}\Phi_{i_2}^{(3)}\Phi_{i_3}^{(3)}p_{i_1+i_2+i_3}+
\tsonethr\Phi_{i_1+i_2+i_3}^{(3)} p_{i_1}p_{i_2}p_{i_3}
+\Phi_{i_1}^{(3)}\Phi_{i_2+i_3}^{(3)}  p_{i_1+i_2}p_{i_3}\right)
=
\frac{\partial \Phi^{(3)}}{\partial u}.
\end{equation}
We do not know of any method for solving this equation for $\Phi^{(3)}$ explicitly,
and have not been able to conjecture the solution from numerical
computations, as we could for $k=2.$ 
However, as we show in the next section, we are able to determine the low
order terms of $\Phi^{(k)}$ in the $p$'s, for arbitrary $k$.


\section{Restriction of the differential equation by grading}

\def\drF#1#2{F^{(#1)}_{i_{#2}}}

In this section we determine
a partial differential equation for
$P_k^{(m)}$ that can be used recursively to construct $P_k^{(m)}$ for
all $m\ge1$. Our method is to apply the symmetrisation operator $\psi_m$ to
the partial differential equation~(\ref{phidiffnew})
given in Theorem~\ref{Tdiff1}.
Some notation is needed for this purpose.
Let
$$h^+_i(x_1,\ldots,x_k)=
\sum_{\atp{j_1,\ldots, j_k\ge1}{j_1+\cdots+j_k=i}}x_1^{j_1}\cdots x_k^{j_k},$$
for $i\geq 1$, and $h^+_0 (x_1,\ldots,x_k)=1$.
Now let 
$ H^+(t;x_1,\ldots,x_k)=\sum_{i\ge0} h^+_i(x_1,\ldots,x_k)t^i.$
If $a(t)$ and $b(t)$ are the generating series for $\{a_i\}$ and $\{b_i\},$
let $a\circ b$ denote the summation $\sum_{i\geq 0}a_i b_i$. This is essentially
the {\em umbral} composition of $a$ and $b$ with respect to $t$,
which will be the only indeterminate used in this paper for umbral composition.

We begin by showing in two particular examples how the action of $\psi_m$,
defined in~(\ref{psidef}),
can be expressed conveniently in terms of umbral composition.
These will suffice to indicate the general procedure.
The first example is the application of $\psi_2$ to the
partial differential equation~(\ref{twocase}). As a preliminary,
apply $\psi_2$ to the final term on the left hand side, which yields
$$\sum_{i_1,i_2\ge1}
x_1^{i_1}x_2^{i_2}(i_1+i_2)\frac{\partial}{\partial p_{i_1+i_2}}F_2^{(1)}
=\sum_{j\geq 1}\left( \sum_{\atp{i_1,i_2\ge1}{i_1+i_2=j} }x_1^{i_1}x_2^{i_2}\right)
j\frac{\partial}{\partial p_j}F_2^{(1)}=
t \frac{\partial}{\partial t} P^{(1)}_2(t)\circ H^+(t;x_1,x_2).$$
Note that the presence
of umbral composition in this expression is explained in an entirely
elementary way.
The other terms require no explanation, and application of $\psi_2$ to the
partial differential equation~(\ref{twocase}) yields
$$\tsh x_1\frac{\partial}{\partial x_1}P^{(1)}_2(x_1)
x_1\frac{\partial P^{(2)}_2}{\partial x_1}
+\tsh x_2\frac{\partial}{\partial x_2}P^{(1)}_2(x_2)
x_2\frac{\partial P^{(2)}_2}{\partial x_2}
+  t \frac{\partial}{\partial t} P^{(1)}_2(t)\circ H^+(t;x_1,x_2)$$
$$= \left(x_1\frac{\partial}{\partial x_1}
+x_2\frac{\partial}{\partial x_2}\right)
P^{(2)}_2.$$
Now, for more variables and more complicated equations the
main complication is the proliferation of terms that arises
from adding the contributions from permuting the variables. To organize this
we introduce another symmetrisation operator, 
$\Xi_m$, defined on power series in $x_1,\ldots ,x_m$ by
$$\Xi_m f(x_1,\ldots,x_m) = {\sum_{\sigma\in\symgp_m}}^{\!\!\!\prime}
f(x_{\sigma(1)},\ldots,x_{\sigma(m)}),$$ where the ``$\,\mbox{}^\prime\,$''
indicates that the summation
is over distinct terms.


For the second example, illustrating this complication,
we apply $\psi_4$ just to the third term on the left hand side
of the partial differential equation~(\ref{thrpde}), to obtain
$$\Xi_4\left(
x_1\frac{\partial}{\partial x_1}P^{(3)}_3(x_1,x_3,x_4)
t \frac{\partial}{\partial t}P^{(1)}_3(t)\circ H^+(t;x_1,x_2)
+x_1\frac{\partial}{\partial x_1}P^{(2)}_3(x_1,x_3)
t \frac{\partial}{\partial t}P^{(2)}_3(t,x_4)\circ H^+(t;x_1,x_2)\right.$$
$$\left. +x_1\frac{\partial}{\partial x_1}P^{(1)}_3(x_1)
t \frac{\partial}{\partial t}P^{(3)}_3(t,x_3,x_4)\circ H^+(t;x_1,x_2)
\right).$$

Without further discussion, in the following result, we now
apply the operator $\psi_m$ to the partial differential
equation~(\ref{phidiffnew})
given in Theorem~\ref{Tdiff1} directly, 
yielding a partial differential equation for 
$P_k^{(m)}$, for
all $m\ge1$. This equation is in terms of $P_k^{(1)},\ldots ,P_k^{(m-1)}$, 
for $m\geq 2$, and we will use it recursively, starting with $P_k^{(1)}$,
in Section 4.
In the statement of the result, let
$$\Psi(z,x_1;)=\sum_{m\ge1}x_1\frac{\partial}{\partial x_1}P^{(m)}(x_1;)z^{m-1}.$$
In addition, we adopt the convention that
$P^{(i)}(x_1;)P^{(j)}(x_1;)$ denotes
$P^{(i)}(x_1,\alpha)P^{(j)}(x_1,\beta)$, for $i+j=n$,
where $(\alpha,\beta)$ is a canonical bipartition of
$\{1,\ldots,n\}-\{1,2\}$ of size $(i-1,j-1).$
Let $\cN_{\widehat{\frkt}}(w)$ be the neighbour set of $w\in\cV_{\widehat{\frkt}}.$


\begin{theorem}\label{TPm}
For $m\geq 1$, the partial differential equation for $P^{(m)}$ is
\begin{eqnarray*}
&\mbox{}&
\left[y^kz^m\right]\,
\sum_{\widehat{\frkt}}\frac{\Xi_m}{\aut(\widehat{\frkt})}
\left(\prod_{\sv\in\cV_{black}(\widehat{\frkt})}\!\!\!\!\!\!
z\left(\frac{y}{1-y\Psi(z,x_\sv;)}\right)^{d_{\widehat{\frkt}}(\sv)}\right) 
\left(\prod_{\sw\in\cV_{white}(\widehat{\frkt})}\!\!\!\!\!\!
\Psi(z,t;)\circ H^+(t;\cN_{\widehat{\frkt}}(\sw))\right) \\
&\mbox{}&=\frac{1}{k-1} \left(
x_1\frac{\partial}{\partial x_1}+\cdots+x_m\frac{\partial}{\partial x_m}
+m-2\right)P^{(m)}.
\end{eqnarray*}
\end{theorem}
\proof  Let $m$ denote the number of black vertices in $\widehat{\frkt}.$
Assign the symbols $1,\ldots,m$ arbitrarily to these vertices.

With each black vertex $\sv$ of $\widehat{\frkt}$ associate the expression
$$
z\left(\frac{y}{1-y\Psi(z,x_\sv;)}\right)^{d_{\widehat{\frkt}}(\sv)}.
$$
This accounts for the attachment of monovalent white vertices by
edges to $\sv.$

With each white vertex $\sw$ of $\widehat{\frkt}$ associate the expression
$$
\Psi(z,t;)\circ H^+(t;\cN_{\widehat{\frkt}}(\sw)).
$$
The result follows
by taking the product of these expressions. \qed

\medskip
The application of the coefficient operator $[y^kz^m]$ is routine
but increasingly laborious as $m$ increases. In the next section we
will carry this out for $m=1,2,3$. The following
result will be needed to give explicit forms for the umbral composition
with $H^+$.



\begin{proposition}\label{CfH}
Let $f(t)$ be a formal power series in $t.$
Then
$$
f(t)\circ H^+(t;x_1,\ldots,x_k)=
\sum_{i=1}^k
f(x_i){\prod_{\atp{1<p\le k}{p\neq i}}}\frac{x_p}{x_i-x_p}.
$$
\end{proposition}


\section{Proofs of the supporting theorems}

\subsection{Proof of Theorem~\ref{Tonept}}
Consider the case $m=1$ in
Theorem~\ref{TPm}. Then since contributions on
the left hand side come only from the tree
$\widehat{\frkt_1},$  we obtain the differential equation
$$
\frac{1}{k} \left( x_1\frac{d P^{(1)}}{d x_1} \right)^k
= \frac{1}{k-1}\left( x_1\frac{d}{d x_1} -1\right) P^{(1)}
$$
for $P^{(1)}.$
To solve this equation, differentiate the equation with respect to $x_1$ and
multiply by $x_1$. Then, with
$f= x_1{d P^{(1)}}/{d x_1}$, we obtain
$$f^{k-1}x_1\frac{d f}{d x_1}
=\frac{1}{k-1}x_1\frac{d f}{d x_1}
-\frac{1}{k-1}f,$$
so, solving for $x_1\frac{d f}{d x_1}$, we have
$$x_1\frac{d f}{d x_1}=\frac{f}{1-(k-1)f^{k-1}}.$$
It is now straightforward to determine, for formal power series in $x$, that $f=w_1$,
by comparing this differential equation with~(\ref{equAA}), and
using the initial condition $f(0)=0$.
The result follows immediately. \qed


\subsection{Proof of Theorem~\ref{Ttwopt}}
Consider  the case $m=2$ in
Theorem~\ref{TPm}. Now contributions on
the left hand side come only from the trees
$\widehat{\frkt_1}$ and $\widehat{\frkt_2}$.
Thus, 
substituting the expression for $P^{(1)}$ from Theorem~\ref{Tonept},
and applying Proposition~\ref{CfH} to carry out the umbral compositions,
we obtain
\begin{eqnarray*}
w_1^{k-1} x_1\frac{\partial P^{(2)}}{\partial x_1}
+ w_2^{k-1} x_2\frac{\partial P^{(2)}}{\partial x_2}
+\frac{x_2w_1-x_1 w_2}{x_1-x_2}\,
\frac{w_1^{k-1}-w^{k-1}_2}{w_1-w_2}
=\frac{1}{k-1} \left(x_1\frac{\partial}{\partial x_1}+x_2\frac{\partial}{\partial x_2}\right)
P^{(2)}.
\end{eqnarray*}
so, rearranging, we have
\begin{eqnarray*}
\frac{1}{k-1}\left( (1-(k-1)w_1^{k-1}) x_1\frac{\partial }{\partial x_1}
+(1-(k-1)w_2^{k-1}) x_2\frac{\partial }{\partial x_2}\right) P^{(2)}
=\frac{x_2w_1-x_1w_2}{x_1-x_2}\,
\frac{w_1^{k-1}-w_2^{k-1}}{w_1-w_2}.
\end{eqnarray*}
\def\wopk
{\left(w_1\frac{\partial}{\partial w_1}+w_2\frac{\partial}{\partial w_2}\right)}
It is now straightforward to verify that
$$P^{(2)} (x_1,x_2)=\log\left( \frac{w_1-w_2}{x_1-x_2} \right)-
\frac{w_1^k - w_2^k}{w_1 - w_2},$$
by confirming that it satisfies the above differential equation, and the
initial condition $P^{(2)} (0,0)=0$. (Note that the constant term in the expansion
of $(w_1-w_2)/(x_1-x_2)$ as a formal power series in $x_1,x_2$ is 1, so 
the logarithm exists.)

Finally, apply the operator $x_1\frac{\partial}{\partial x_1}+
x_2\frac{\partial}{\partial x_2}$ to $P^{(2)}$,
and the result follows.  \qed


\subsection{Proof of Theorem~\ref{Tthreept}}
Consider  the case $m=3$ in
Theorem~\ref{TPm}. Then since contributions on
the left hand side come only from the trees
$\widehat{\frkt_1},\ldots,\widehat{\frkt_4},$
having substituted the expression for $P^{(1)}$ from Theorem~\ref{Tonept},
it follows that

\begin{eqnarray*}
&\mbox{}& \Xi_3 w_1^{k-1} x_1\frac{\partial}{\partial x_1} P^{(3)} \\
&\mbox{}&+\Xi_3 (k-1) w_1^{k-2} \left(x_1\frac{\partial}{\partial x_1} P^{(2)}(x_1,x_2)\right)\,
\left(x_1\frac{\partial}{\partial x_1} P^{(2)}(x_1,x_3)\right)  \\
&\mbox{}&+\Xi_3 \left(\frac{\partial}{\partial w_1}
\frac{w_1^{k-1}-w_2^{k-1}}{w_1-w_2}\right)\,
\left(x_1\frac{\partial}{\partial x_1} P^{(2)}(x_1,x_3)\right)\,
\left( w(t)\circ H^+(t;x_1,x_2)\right) \\
&\mbox{}&+\Xi_3
\left(\frac{w_1^{k-1}-w_2^{k-1}}{w_1-w_2}\right)\,
\left(t\frac{\partial}{\partial t} P^{(2)}(t,x_3)\circ H^+(t;x_1,x_2) \right)\\
&\mbox{}&+\Xi_3 2h_{k-3}(w_1,w_2,w_3)
\left(w(t)\circ  H^+(t;x_1,x_2,x_3)\right) \\
&\mbox{}&+\Xi_3 \left(\frac{\partial}{\partial w_2}
h_{k-3}(w_1,w_2,w_3)\right)\,
\left(w(t)\circ  H^+(t;x_1,x_2)\right) \,
\left(w(t)\circ  H^+(t;x_2,x_3)\right)  \\
&=& \frac{1}{k-1}\left(
x_1\frac{\partial}{\partial x_1}+x_2\frac{\partial}{\partial x_2}+x_3\frac{\partial}{\partial x_3}
+1\right) P^{(3)}.
\end{eqnarray*}
The six expressions on the left hand side arise  from $\widehat{\frkt_1},$
$\widehat{\frkt_1},$ $\widehat{\frkt_2},$ $\widehat{\frkt_2},$
 $\widehat{\frkt_3}$ and $\widehat{\frkt_4},$
respectively. 
Note that, under the action of $\Xi_3,$ the six expressions
on the left hand side expand into $3, 3, 6, 3, 1$ and $3$ terms, respectively.
Now apply Proposition~\ref{CfH} to carry out the umbral compositions, and
use the fact that
\begin{equation}\label{xtow}
(1-(k-1)w^{k-1})x\frac{\partial}{\partial x}=
w\frac{\partial}{\partial w}
\end{equation}
(this latter follows from~(\ref{equAA})).
Simplifying with the help of Maple, we obtain
\def\ww#1#2#3{\frac{w_#1^{k-1}-w_#2^{k-1}}{(w_#1-w_#2)^2}\,(w_#2A_{#1#3}-w_#1A_{#2#3})}
\begin{eqnarray*}
&\mbox{}& \frac{1}{k-1}\left(\sum_{i=1}^3
w_i\frac{\partial}{\partial w_i}+1\right) P^{(3)} =
(k-1)\left(w_1^{k-2}A_{12}A_{13} + w_2^{k-2}A_{21}A_{23} +w_3^{k-2}A_{31}A_{32}\right) \\
&\mbox{}& + \ww123 + \ww132\\ 
&\mbox{}&+ \ww231,
\end{eqnarray*}
where
$$
A_{ij}=\frac{w_iw_j}{1-(k-1)w_i^{k-1}}\,
\frac{w_i^{k-1}-w_j^{k-1}}{(w_i-w_j)^2}.
$$
The solution to this equation is given in Theorem~\ref{Tthreept}, and has
been verified with the aid of Maple, giving the desired result. \qed


\section{Computational comments and conjectures}
We have shown in Section 4 that $P^{(1)}, P^{(2)}$ and $P^{(3)}$ can each be obtained
as the solutions to first order linear partial differential equations.
We believe that $P^{(m)},$ for $m\ge4,$ can be obtained in a similar way as 
the solution of  such an equation. Moreover, we conjecture that the equation
for any $m\geq 3$,(obtained from Theorem~\ref{TPm}, and applying~(\ref{xtow})
as described for $m=1,2,3$ in Section 4)
after multiplying through by $k-1$, is of the form
$$\left(\sum_{i=1}^m
w_i\frac{\partial}{\partial w_i}+(m-2)\right) P^{(m)}
=R_m(w_1,\ldots ,w_m),$$
where $R_m$ is a rational function in $w_1,\ldots ,w_m$, obtained
from $P^{(1)},\ldots, P^{(m-1)}$. That is, there is no
dependency of $R_m$ on $x_1 ,\ldots ,x_m$ except through~(\ref{Ewx}). Now let
$Q^{(m)}(t)$ be obtained by substituting $tw_i$ for $w_i$ in
$P^{(m)}$  for $1=1,\ldots ,m$. Then the above partial differential
equation is transformed into the first order linear ordinary differential
equation
\begin{equation}\label{tdiffeq}
{\frac{d}{dt}}(t^{m-2}Q^{(m)}(t))=t^{m-3}R_m(tw_1,\ldots ,tw_m),
\end{equation}
which can be solved routinely, in theory. In practice, this is precisely
how we obtained $P^{(3)}$, with the aid of Maple, in Section 4 above. However,
even in this case, the simplification of the equation was difficult; we
provided human help by proving that the rational expression on the right hand
side of the equation is independent of the $x$'s, and then
replaced each $x_i$ by $w_i$ to evaluate it. This explains how the
$A_{ij}$ arise, 
as $x_i\frac{\partial} {\partial x_i} P^{(2)}(x_i,x_j)$
evaluated at $x_i=w_i$ and $x_j=w_j$.

For $m=4$,
the expressions became too big to be tractable, and
we have not found a convenient way of circumventing this. We conjecture
that, for each $m\geq 3$, $P^{(m)}$ is a rational function of
$w_1,\ldots ,w_m$, whose denominator is consistent with
Conjecture~\ref{CConj}, using~(\ref{equAA}).
(Note that for $m=2$, the right hand side of the equation, as
obtained in the Proof of Theorem~\ref{Ttwopt}, is not a
rational function of $w_1,w_2$ alone, but rather involves $x_1,x_2$
also.)
%
%

\section*{Acknowledgements}
This work was supported by grants from the Natural
Sciences and Engineering Research Council of Canada.


\end{document}